 \documentclass[final,3p,times]{elsarticle}
\usepackage{amssymb}
 \usepackage{amsthm}
 \usepackage[T1]{fontenc}
 
\usepackage[ansinew]{inputenc} 
\usepackage{mathrsfs}
\usepackage{euscript}
\usepackage{natbib}

\newtheorem{theorem}{Theorem}%[section]
\newtheorem{lemma}{Lemma}%[section]
\newtheorem{proposition}{Proposition}%[section]
\newdefinition{rmk}{ Remark}%[section]
\newproof{pf}{Proof}

\begin{document}

\begin{frontmatter}


%\tnotetext[label1]{}
%% \ead{ahmad.fino01@gmail.com}

\title{ Finite time blow up for wave equations with strong damping in an exterior domain}

 \author[label]{Ahmad Z. FINO}
 \address[label]{LaMA-Liban, Lebanese University, Faculty of Sciences, Department of Mathematics, P.O. Box 826 Tripoli, Lebanon}

%\author{}

%\address{}
 \ead{ahmad.fino01@gmail.com; afino@ul.edu.lb}
 \begin{abstract}
We consider the initial boundary value problem in exterior domain for semilinear wave equations with power-type nonlinearity $|u|^p$. We will establish blow-up results when $p$ is less than or equal to Strauss' exponent which is the same one for the whole space case $\mathbb{R}^n$.
\end{abstract}

\begin{keyword}
Semilinear wave equation\sep Blow-up \sep Initial-boundary value problem\sep Exterior domain\sep Strong damping

\MSC[2010] 35L05 \sep 35L70 \sep 35B33 \sep 34B44
\end{keyword}

\end{frontmatter}

% \linenumbers

\section{Introduction}
\setcounter{equation}{0} 

This paper concerns the initial boundary value problem of the strongly damped wave equation in an exterior domain. Let $\Omega\subset\mathbb{R}^n$ be an exterior domain whose obstacle 
$\mathcal{O}\subset\mathbb{R}^n$ is bounded with smooth compact boundary $\partial\Omega$. We consider the initial boundary value problem
\begin{equation}\label{eq1}
\left\{\begin{array}{ll}
\,\, \displaystyle {u_{tt}-\Delta
u -\Delta u_t =|u|^p} &\displaystyle {t>0,x\in {\Omega},}\\
{}\\
\displaystyle{u(0,x)= \varepsilon u_0(x),\;\;u_t(0,x)= \varepsilon u_1(x)\qquad\qquad}&\displaystyle{x\in {\Omega},}\\
{}\\
\displaystyle{u=0,\qquad\qquad}&\displaystyle{t\geq0,\;x\in {\partial\Omega},}
\end{array}
\right.
\end{equation} 
where the unknown function $u$ is real-valued, $n\geq 1$, $\varepsilon>0$, and $p>1.$ Throughout this paper, we assume that 
\begin{equation}\label{condition1}
(u_0,u_1)\in H^1_0(\Omega)\times L^2(\Omega),\quad\mbox{and}\quad u_0,u_1\geq0.
\end{equation}
Without loss of generality, we assume that $0\in\mathcal{O}\subset\subset B(R)$, where $B(R):=\{x\in\mathbb{R}^n:\;|x|< R\}$ is a ball of radius $R$ centred at the origin, and that  
\begin{equation}\label{condition2}
\mbox{supp}u_i\subset B(R),\quad i=0,1.
\end{equation}
For the simplicity of notations, $\|\cdotp\|_q$ and $\|\cdotp\|_{H^1}$ $(1\leq q\leq \infty)$ stand for the usual $L^q(\Omega)$-norm and $H^1_0(\Omega)$-norm, respectively.

First, the following local well-posedness result is needed. 
\begin{proposition}\label{prop1}\cite[see Proposition 2.1]{IkehataInoue}\\
Let $1<p<\infty$. Under the assumptions $(\ref{condition1})$-$(\ref{condition2})$, there exists a maximal existence time $T_{\max}>0$ such that the problem $(\ref{eq1})$ possesses a unique weak solution
$$u\in C([0,T_{\max}),H^1_0(\Omega))\cap C^1([0,T_{\max}), L^2(\Omega)),$$
where $0< T_{\max}\leq\infty.$ Moreover, $u(t,\cdotp)$ is supported in the ball $B(t+R).$ In addition:
\begin{equation}\label{alternative}
\mbox{either}\;\;T_{\max}=\infty \quad\mbox{or else}\quad T_{\max}<\infty \;\;\mbox{and}\;\; \|u(t,\cdotp)\|_{H^1_0}+\|u_t(t,\cdotp)\|_2\rightarrow\infty\;\;\mbox{as}\;\; t\rightarrow T_{\max}.
 \end{equation}
\end{proposition}
\begin{rmk}
We say that $u$ is a global solution of $(\ref{eq1})$ if $T_{\max}=\infty,$ while in the case of $T_{\max}<\infty,$ we say that $u$ blows up in finite time.
\end{rmk}
\noindent Let $p_c(n)=+\infty$, for $n=1$, and let $p_c(n)$, for $n\geq 2$, be the positive root of the quadratic equation
$$(n-1)p^2-(n+1)p-2=0.$$
The number $p_c(n)$ is known as the critical exponent (Strauss exponent) of the semilinear wave equation
\begin{equation}\label{}
\left\{\begin{array}{ll}
\,\, \displaystyle {u_{tt}-\Delta
u=|u|^p} &\displaystyle {t>0,x\in {\Omega},}\\
{}\\
\displaystyle{u(0,x)= u_0(x),\;\;u_t(0,x)= u_1(x)\qquad\qquad}&\displaystyle{x\in {\Omega},}\\
{}\\
\displaystyle{u=0,\qquad\qquad}&\displaystyle{t>0,\;x\in {\partial\Omega},}
\end{array}
\right.
\end{equation} 
since it divides $(1,\infty)$ into two subintervals such that the following description holds: If $p\in (1, p_c(n))$, then solutions with nonnegative initial values blow-up in finite time; if $p\in(p_c (n),\infty)$, then solutions with small (and sufficiently regular) initial values exist for all time (see e.g. \cite{Strauss1}). The proof
has an interesting and exciting history that spans three decades. We only give a brief summary here and refer the reader to \cite{Strauss1,Jiao} and the references therein for details. The problem as regards the existence or nonexistence of global solutions is sometimes referred to as the Conjecture of Strauss \cite{Strauss2}. The same problem was also posed by Glassey \cite{Glassey}.

Our main result is the following
\begin{theorem}\label{theo1}
 Assume that the initial data satisfy $(\ref{condition1})$-$(\ref{condition2})$. If 
 $$\left\{\begin{array}{l}
 1<p< p_c(1)=\infty,\\
  1<p\leq p_c(2),\\
   1<p< p_c(n),\quad n\geq3,
 \end{array}
 \right.
 $$
 then the solution of the problem $(\ref{eq1})$ blows up in finite time.
\end{theorem}
\rmk It still an open problem to prove that the solution of the problem $(\ref{eq1})$ blows up in finite time for $p=p_c(n)$, $n\geq3$.\\

This paper is organized as follows: in Section \ref{sec2}, we present several preliminaries. Section \ref{sec3} contains the proofs of the blow-up theorem (Theorem \ref{theo1}).

%%%%%%%%%%%%%%%%%%%%%%%%%%%%%%%%%%%%     Section 2     %%%%%%%%%%%%%%%%%%%%%%%%%%%%%%%%%%%%%

\section{Preliminaries}\label{sec2}
\setcounter{equation}{0}
In this section, we give some preliminary properties that will be used in the proof of Theorem \ref{theo1}.\\
In \cite[p.~386]{Sideris}, Sideris obtained the following well-known ODE blow-up result:
\begin{lemma}\label{lemma1}
Let $p>1$, $a\geq1$, and $(p-1)a>q-2$. If $F\in C^2([0,T))$ satisfies
\begin{enumerate}
\item $F(t)\geq \delta(t+R)^a$, and
\item $\frac{d^2 F(t)}{dt^2}\geq k(t+R)^{-q}[F(t)]^p$,
\end{enumerate}
with some positive constants $\delta,k,$ and $R$, then $T<\infty$.
\end{lemma}
To prove the main results in this paper when $n=2$, we will concentrate on the
improvement of the above well-known Sideris ODE blow-up result, for when the differential inequality involves a
logarithmic term.
\begin{lemma}\label{lemma10}\cite[Lemma 2.3]{Han}
Let $p>1$, $a\geq1$, and $(p-1)a>q-2$. If $F\in C^2([0,T))$ satisfies
\begin{enumerate}
\item $F(t)\geq \delta(t+R)^a$, and
\item $\frac{d^2 F(t)}{dt^2}\geq k[\ln(t+R)]^{-q/2}(t+R)^{-q}[F(t)]^p$,
\end{enumerate}
with some positive constants $\delta,k,$ and $R$, then $T<\infty$.
\end{lemma}
\begin{lemma}\label{lemma11}\cite[Lemma 2.4]{Han}
Let $p>1$, $a\geq1$, and $(p-1)a=q-2$. If $F\in C^2([0,T))$ satisfies, when $t\geq T_0>0$,
\begin{enumerate}
\item $F(t)\geq K_0(t+R)^a$, and
\item $\frac{d^2 F(t)}{dt^2}\geq K_1[\ln(t+R)]^{-q/2}(t+R)^{-q}[F(t)]^p$,
\end{enumerate}
with some positive constants $K_0,K_1,T_0$ and $R$. Fixing $K_1$, there exists a positive constant $c_0$, independent of $R$ and $T_0$, such that if $K_0\geq c_0$, then $T<\infty$.
\end{lemma}
We also need of the following special functions.
\begin{lemma}\label{lemma3}\cite[Lemma 2.2]{ZhouHan}
There exists a function $\phi_0(x)\in C^2(\Omega)$ satisfying the following boundary value problem
\begin{equation}\label{function1}
\left\{\begin{array}{l}
\Delta \phi_0(x)=0, \;\mbox{in}\;\Omega,\;\;n\geq 3,\\
\phi_0|_{\partial\Omega}=0,\\
|x|\rightarrow\infty,\quad \phi_0(x)\rightarrow 1.
\end{array}
\right.
\end{equation}
Moreover, $\phi_0(x)$ satisfies: for all $x\in\Omega$, $0<\phi_0(x)<1$.
\end{lemma}
\begin{lemma}\label{lemma4}\cite[Lemma 2.5]{Han}
There exists a function $\phi_0(x)\in C^2(\Omega)$ satisfying the following boundary value problem
\begin{equation}\label{function2}
\left\{\begin{array}{l}
\Delta \phi_0(x)=0, \;\mbox{in}\;\Omega,\;\;n=2,\\
\phi_0|_{\partial\Omega}=0,\\
|x|\rightarrow\infty,\quad \phi_0(x)\rightarrow +\infty,\;\;\mbox{and}\;\phi_0(x)\;\mbox{increase at the rate of}\;\ln (|x|).
\end{array}
\right.
\end{equation}
Moreover, $\phi_0(x)$ satisfies: for all $x\in\Omega$, $0<\phi_0(x)\leq C\ln r$, where $r=|x|$ and $C>0$ is a positive contant.
\end{lemma}
\begin{lemma}\label{lemma5}\cite[Lemma 2.2]{Han2}
There exists a function $\phi_0(x)\in C^2([0,\infty))$ satisfying the following boundary value problem
\begin{equation}\label{function2}
\left\{\begin{array}{l}
\Delta \phi_0(x)=0, \;x>0,\\
\phi_0|_{x=0}=0,\\
x\rightarrow\infty,\quad \phi_0(x)\rightarrow +\infty,\;\;\mbox{and}\;\phi_0(x)\;\mbox{increase at the rate of linear function $x$}.
\end{array}
\right.
\end{equation}
Moreover, $\phi_0(x)$ satisfies: there exist two positive constants $C_1$ and $C_2$ such that, for all $x>0$, we have $C_1x\leq \phi_0(x)\leq C_2x$.
\end{lemma}
Similarly, we have the following
\begin{lemma}\label{lemma6}
There exists a function $\varphi_1(x)\in C^2(\Omega)$ satisfying the following boundary value problem
\begin{equation}\label{function3}
\left\{\begin{array}{l}
\Delta \varphi_1(x)=\frac{1}{2}\varphi_1(x), \;\mbox{in}\;\Omega,\;\;n\geq 1,\\
\varphi_1|_{\partial\Omega}=0.
\end{array}
\right.
\end{equation}
Moreover, $\varphi_1(x)$ satisfies: there exists positive constant $C_1$, for all $x\in\Omega$, $0<\varphi_1(x)\leq C_1(1+|x|)^{-(n-1)/2}e^{|x|}$.
\end{lemma}
\proof It is sufficient to take $\varphi_1(x)=\phi_1(\frac{x}{\sqrt{2}})$ where $\phi_1$ is the function defined by \cite[Lemma 2.3]{ZhouHan} on $\frac{1}{\sqrt{2}}\Omega$ instead of $\Omega$. \hfill $\square$\\

In order to continue the description of the following lemmas, we define the following test function
$$\psi_1(x,t)=\varphi_1(x)e^{-t},\quad \forall \;x\in\Omega,\;t\geq0.$$
It is easy to check that
$$
(\psi_1)_t(x,t)=-\psi_1(x,t),\qquad (\psi_1)_{tt}(x,t)=\psi_1(x,t),\qquad \mbox{and}\quad \Delta\psi_1(x,t)=\frac{1}{2}\psi_1(x,t).
$$
\begin{lemma}\label{lemma7}\cite[Lemma 2.4]{ZhouHan}\\
Let $p>1$, $n\geq1$. Then, for all $t\geq0$, we have
$$\int_{\Omega\cap\{|x|\leq t+R\}}[\psi_1(x,t)]^{p'}\,dx\leq C(t+R)^{n-1-(n-1)p'/2},$$
where $p'=p/(p-1)$ and $C$ is a positive constant.
\end{lemma}
\begin{lemma}\label{lemma8}
Let $p>1$, $n\geq1$. Then, for all $t\geq0$, we have
$$\int_{\Omega\cap\{|x|\leq t+R\}}[\phi_0(x)]^{-1/(p-1)}[\psi_1(x,t)]^{p'}\,dx\leq C(t+R)^{n-1-(n-1)p'/2},$$
where $p'=p/(p-1)$ and $C$ is a positive constant.\\
For the case $n=2$, we can improve the last inequality, more precisely, there exists $R_1\gg1$ such that, for all $t\geq R_1$, we have
$$\int_{\Omega\cap\{|x|\leq t+R\}}[\phi_0(x)]^{-1/(p-1)}[\psi_1(x,t)]^{p'}\,dx\leq C(t+R)^{1-p'/2}(\ln(t+R))^{-1/(p-1)}.$$
\end{lemma}
\proof For the case $n\geq 3$, see \cite[Lemma 2.5]{ZhouHan}. For the case of $n=2$ see \cite[Lemma 2.8]{Han}. Finally, for the one dimensional case see \cite[Lemma 2.5]{Han2}.\hfill$\square$
%%%%%%%%%%%%%%%%%%%%%%%%%%%%%%%%%%%%%%       Section 3       %%%%%%%%%%%%%%%%%%%%%%%%%%%%%%%%%%%%%%

\section{Proof of Theorem \ref{theo1}}\label{sec3}
\setcounter{equation}{0}
Theorem \ref{theo1} is a consequence of the lower bound and the blowup result about nonlinear differential inequalities in Lemmas \ref{lemma1} and \ref{lemma2}.

To outline the method, we will introduce the following functions:
$$
\left\{\begin{array}{l}
\displaystyle F_0(t)=\int_{\Omega}u(x,t)\phi_0(x)\,dx,\\
{}\\
\displaystyle F_1(t)=\int_{\Omega}u(x,t)\psi_1(x,t)\,dx.
\end{array}
\right.
$$
By density we can assume that the solution $u$ is sufficiently smooth, which implies that $F_0(t)$ and $F_1(t)$ are well-defined $C^2$-functions for all $t\geq0$. The following lemma is dedicated to obtain a lower bound on $F_1(t)$.
\begin{lemma}\label{lemma9}
Let $n\geq1$. Under the assumptions $(\ref{condition1})$-$(\ref{condition2})$, let $(u,u_t)$ be the solution of the problem $(\ref{eq1})$ such that
$$(u,u_t)\in C([0,T),H^1_0(\Omega))\times C([0,T), L^2(\Omega)),$$
and
$$ \mbox{supp}(u,u_t)\subset B(t+R):=\{x\in\Omega:\;|x|< t+R\}.$$
Then, for all $t\geq0$, we have
$$F_1(t)\geq \left( \frac{\varepsilon}{3}\left(1-e^{-\frac{3}{2}t}\right)+\varepsilon e^{-\frac{3}{2}t}\right)\int_{\Omega}\varphi_1(x)u_0(x)\,dx+ \frac{2\varepsilon}{3}\left(1-e^{-\frac{3}{2}t}\right)\int_{\Omega}\varphi_1(x)u_1(x)\,dx\geq \varepsilon c_0>0.$$
\end{lemma}
\proof We multiply $(\ref{eq1})$ by the test function $\psi_1\in C^2(\Omega\times\mathbb{R})$ and integrate over $\Omega\times[0,t]$, we get
\begin{equation}\label{equation1}
\int_0^t\int_{\Omega}u_{ss}\psi_1\,dx\,ds-\int_0^t\int_{\Omega}\Delta u\psi_1\,dx\,ds-\int_0^t\int_{\Omega}\Delta u_s\psi_1\,dx\,ds=\int_0^t\int_{\Omega}|u|^p\psi_1\,dx\,ds.
\end{equation}
Use integration by parts, we have
\begin{eqnarray*}
\int_0^t\int_{\Omega}u_{ss}\psi_1\,dx\,ds&=&-\int_0^t\int_{\Omega}u_{s}(\psi_1)_s\,dx\,ds+\int_{\Omega}u_{t}(x,t)\psi_1(x,t)\,dx-\varepsilon \int_{\Omega}u_1(x)\varphi_1(x)\,dx\\
&=&\int_0^t\int_{\Omega}u(x,s)\psi_1(x,s)\,dx\,ds+\int_{\Omega}[u(x,t)+u_t(x,t)]\psi_1(x,t)\,dx-\varepsilon\int_{\Omega}[u_{0}(x)+u_1(x)]\varphi_1(x)\,dx,\\
\end{eqnarray*}
where we have used the fact that $(\psi_1)_s(x,s)=-\psi_1(x,s)$ and $(\psi_1)_{ss}(x,s)=\psi_1(x,s)$, for all $x\in\Omega$, $s\geq0$. Moreover
\begin{eqnarray*}
\int_0^t\int_{\Omega}\Delta u\psi_1\,dx\,ds&=&-\int_0^t\int_{\Omega}\nabla u\nabla \psi_1\,dx\,ds+\int_0^t\int_{\partial\Omega}\psi_1\nabla u\cdotp\bold{n}\,d\sigma\,ds\\
&=&\int_0^t\int_{\Omega} u\Delta \psi_1\,dx\,ds-\int_0^t\int_{\partial\Omega}u\nabla \psi_1\cdotp\bold{n}\,d\sigma\,ds+\int_0^t\int_{\partial\Omega}\psi_1\nabla u\cdotp\bold{n}\,d\sigma\,ds,\\
&=&\frac{1}{2}\int_0^t\int_{\Omega} u \psi_1\,dx\,ds,
\end{eqnarray*}
where we have used the boundary conditions and the fact that $\Delta \psi_1(x,s)=\frac{1}{2}\psi_1(x,s)$, for all $x\in\Omega$, $s\geq0$. Similarly
\begin{eqnarray*}
\int_0^t\int_{\Omega}\Delta u_{s}\psi_1\,dx\,ds&=&-\int_0^t\int_{\Omega}\Delta u(\psi_1)_s\,dx\,ds+\int_{\Omega}\psi_1(x,t)\Delta u(x,t)\,dx- \int_{\Omega}\varphi_1(x)\Delta u(x,0)\,dx\\
&=&\int_0^t\int_{\Omega}\Delta u \psi_1\,dx\,ds+\int_{\Omega}\psi_1(x,t)\Delta u(x,t)\,dx- \int_{\Omega}\varphi_1(x)\Delta u(x,0)\,dx\\
&=&\frac{1}{2}\int_0^t\int_{\Omega} u \psi_1\,dx\,ds+\frac{1}{2}\int_{\Omega}\psi_1(x,t) u(x,t)\,dx- \frac{1}{2}\varepsilon\int_{\Omega}\varphi_1(x)u_0(x)\,dx,
\end{eqnarray*}
where a similar calculation as above was applied. Combining the above equalities, we conclude from (\ref{equation1}) that
\begin{eqnarray*}
\int_0^t\int_{\Omega}|u|^p\psi_1\,dx\,ds&=&\int_0^t\int_{\Omega}u\psi_1\,dx\,ds+\int_{\Omega}[u(x,t)+u_t(x,t)]\psi_1(x,t)\,dx-\varepsilon\int_{\Omega}[u_{0}(x)+u_1(x)]\varphi_1(x)\,dx\\
&{}&-\frac{1}{2}\int_0^t\int_{\Omega} u \psi_1\,dx\,ds\\
&{}&-\frac{1}{2}\int_0^t\int_{\Omega} u \psi_1\,dx\,ds-\frac{1}{2}\int_{\Omega}\psi_1(x,t) u(x,t)\,dx+ \frac{1}{2}\varepsilon\int_{\Omega}\varphi_1(x)u_0(x)\,dx\\
&=&\frac{d}{dt}\int_{\Omega}u(x,t)\psi_1(x,t)\,dx+\frac{3}{2}\int_{\Omega}u(x,t)\psi_1(x,t)\,dx- \frac{1}{2}\varepsilon\int_{\Omega}\varphi_1(x)u_0(x)\,dx- \varepsilon\int_{\Omega}\varphi_1(x)u_1(x)\,dx\\
&=&\frac{d}{dt}F_1(t)+\frac{3}{2}F_1(t)- \frac{1}{2}\varepsilon\int_{\Omega}\varphi_1(x)u_0(x)\,dx- \varepsilon\int_{\Omega}\varphi_1(x)u_1(x)\,dx.
\end{eqnarray*}
So by $\psi_1>0$, we have
\begin{eqnarray*}
\frac{d}{dt}F_1(t)+\frac{3}{2}F_1(t)&=& \int_0^t\int_{\Omega}|u|^p\psi_1\,dx\,ds+ \frac{1}{2}\varepsilon\int_{\Omega}\varphi_1(x)u_0(x)\,dx+ \varepsilon\int_{\Omega}\varphi_1(x)u_1(x)\,dx\\
&\geq& \frac{1}{2}\varepsilon\int_{\Omega}\varphi_1(x)u_0(x)\,dx+ \varepsilon\int_{\Omega}\varphi_1(x)u_1(x)\,dx.
\end{eqnarray*}
Multiply the above expression by $\displaystyle e^{\frac{3}{2}t}$, we obtain
$$
\frac{d}{dt}(e^{\frac{3}{2}t}F_1(t))\geq \frac{\varepsilon}{2}e^{\frac{3}{2}t}\int_{\Omega}\varphi_1(x)u_0(x)\,dx+ \varepsilon e^{\frac{3}{2}t}\int_{\Omega}\varphi_1(x)u_1(x)\,dx,
$$
and integrating the last differential inequality over $[0,t]$, we get
$$
e^{\frac{3}{2}t}F_1(t)-F_1(0)\geq \frac{\varepsilon}{3}\left(e^{\frac{3}{2}t}-1\right)\int_{\Omega}\varphi_1(x)u_0(x)\,dx+ \frac{2\varepsilon}{3}\left(e^{\frac{3}{2}t}-1\right)\int_{\Omega}\varphi_1(x)u_1(x)\,dx.
$$
As $F_1(0)=\varepsilon\int_{\Omega}\varphi_1(x)u_0(x)\,dx$, we arrive at
$$
F_1(t)\geq \left( \frac{\varepsilon}{3}\left(1-e^{-\frac{3}{2}t}\right)+\varepsilon e^{-\frac{3}{2}t}\right)\int_{\Omega}\varphi_1(x)u_0(x)\,dx+ \frac{2\varepsilon}{3}\left(1-e^{-\frac{3}{2}t}\right)\int_{\Omega}\varphi_1(x)u_1(x)\,dx\geq \varepsilon c_0>0.
$$
 \hfill$\square$
 
Next, in order to apply Lemma \ref{lemma1} on $F_0(t)$, we multiply $(\ref{eq1})$ by $\phi_0$ and integrate over $\Omega$
$$
\int_{\Omega}u_{tt}\phi_0\,dx-\int_{\Omega}\Delta u\phi_0\,dx-\int_{\Omega}\Delta u_t\phi_0\,dx=\int_{\Omega}|u|^p\phi_0\,dx.
$$
By using integration by parts, boundary conditions and Lemma \ref{lemma3}, we can easily check that
$$\int_{\Omega}\Delta u\phi_0\,dx=\int_{\Omega}u\Delta\phi_0\,dx=0,$$
and 
$$\int_{\Omega}\Delta u_t\phi_0\,dx=\int_{\Omega}u_t\Delta\phi_0\,dx=0.$$
Therefore
\begin{equation}\label{equation3}
\frac{d^2}{dt^2}F_0(t)=\int_{\Omega}u_{tt}\phi_0\,dx=\int_{\Omega}|u|^p\phi_0\,dx.
\end{equation}
To estimate the right-hand side of the last equality, we use H\"older's inequality
\begin{eqnarray*}
\left|\int_{\Omega}u(x,t)\phi_0(x)\,dx\right|&=&\left|\int_{\Omega\cap\{|x|\leq t+R\}}u(x,t)[\phi_0(x)]^{1/p}[\phi_0(x)]^{(p-1)/p}\,dx\right|\\
&\leq&\left( \int_{\Omega\cap\{|x|\leq t+R\}}\left|u(x,t)[\phi_0(x)]^{1/p}\right|^p\,dx\right)^{1/p}\left( \int_{\Omega\cap\{|x|\leq t+R\}}\left|[\phi_0(x)]^{(p-1)/p}\right|^{p'}\,dx\right)^{1/{p'}},
\end{eqnarray*}
 where $p'=p/(p-1)$, and then
 \begin{eqnarray*}
|F_0(t)|^p=\left|\int_{\Omega}u(x,t)\phi_0(x)\,dx\right|^p&\leq&\left( \int_{\Omega\cap\{|x|\leq t+R\}}\left|u(x,t)\right|^p\phi_0(x)\,dx\right)\left( \int_{\Omega\cap\{|x|\leq t+R\}}\phi_0(x)\,dx\right)^{p-1}\\
&\leq&\left( \int_{\Omega}\left|u(x,t)\right|^p\phi_0(x)\,dx\right)\left( \int_{\Omega\cap\{|x|\leq t+R\}}\phi_0(x)\,dx\right)^{p-1}.
\end{eqnarray*}
So
\begin{equation}\label{equation5}
\int_{\Omega}\left|u(x,t)\right|^p\phi_0(x)\,dx\geq \frac{|F_0(t)|^p}{\displaystyle \left( \int_{\Omega\cap\{|x|\leq t+R\}}\phi_0(x)\,dx\right)^{p-1}}.
\end{equation}
At this stage, we distinguish the following four cases.\\

\noindent \underline{\textbf{Case $n\geq3$}}: By lemma \ref{lemma3}, we have $0<\phi_0(x)< 1$, then (\ref{equation5}) implies
$$
\int_{\Omega}\left|u(x,t)\right|^p\phi_0(x)\,dx\geq \frac{|F_0(t)|^p}{\displaystyle \left( \int_{\{|x|\leq t+R\}}1\,dx\right)^{p-1}}=|F_0(t)|^p\left[\mbox{Vol}(\bold B^n)\right]^{-(p-1)}(t+R)^{-n(p-1)},$$
where $\bold B^n$ stands for the unit closed ball in $\mathbb{R}^n$. Combining that above inequalities, we infer that
\begin{equation}\label{equation2}
\frac{d^2}{dt^2}F_0(t)\geq k(t+R)^{-n(p-1)}|F_0(t)|^p,
\end{equation}
where $k=\left[\mbox{Vol}(\bold B^n)\right]^{-(p-1)}>0$. So $F_0$ satisfies the second inequality in Lemma \ref{lemma1}. To provide that $F_0$ is also verifies the first  inequality in Lemma \ref{lemma1}, we relate $\frac{d^2}{dt^2}F_0(t)$ to $F_1$ using again H\"older's inequality
\begin{eqnarray*}
\left|\int_{\Omega}u(x,t)\psi_1(x,t)\,dx\right|&=&\left|\int_{\Omega\cap\{|x|\leq t+R\}}u(x,t)[\phi_0(x)]^{1/p}[\phi_0(x)]^{-1/p}\psi_1(x,t)\,dx\right|\\
&\leq&\left( \int_{\Omega\cap\{|x|\leq t+R\}}\left|u(x,t)[\phi_0(x)]^{1/p}\right|^p\,dx\right)^{1/p}\left( \int_{\Omega\cap\{|x|\leq t+R\}}\left|[\phi_0(x)]^{-1/p}\psi_1(x,t)\right|^{p'}\,dx\right)^{1/{p'}}\\
&\leq&\left( \int_{\Omega}\left|u(x,t)\right|^p\phi_0(x)\,dx\right)^{1/p}\left( \int_{\Omega\cap\{|x|\leq t+R\}}[\phi_0(x)]^{-1/(p-1)}[\psi_1(x,t)]^{p'}\,dx\right)^{1/{p'}},
\end{eqnarray*}
then
\begin{eqnarray*}
|F_1(t)|^p=\left|\int_{\Omega}u(x,t)\psi_1(x,t)\,dx\right|^p&\leq&\left( \int_{\Omega}\left|u(x,t)\right|^p\phi_0(x)\,dx\right)\left( \int_{\Omega\cap\{|x|\leq t+R\}}[\phi_0(x)]^{-1/(p-1)}[\psi_1(x,t)]^{p'}\,dx\right)^{p-1}.
\end{eqnarray*}
So, by using Lemmas \ref{lemma8} and \ref{lemma9}, we get
$$\int_{\Omega}\left|u(x,t)\right|^p\phi_0(x)\,dx\geq \frac{|F_1(t)|^p}{\displaystyle\left( \int_{\Omega\cap\{|x|\leq t+R\}}[\phi_0(x)]^{-1/(p-1)}[\psi_1(x,t)]^{p'}\,dx\right)^{p-1}}\geq \frac{(\varepsilon c_0)^p}{\displaystyle \left( C(t+R)^{n-1-(n-1)p'/2}\right)^{p-1}}=L(t+R)^{-(n-1)(p/2-1)},$$
where $L>0$ is a positive contant independent of $t$. Therefore, (\ref{equation3}) implies
$$\frac{d^2}{dt^2}F_0(t)\geq L(t+R)^{-(n-1)(p/2-1)}.$$
Integrate twice, we have
$$F_0(t)\geq \delta(t+R)^{n+1-(n-1)p/2}+\frac{dF_0(0)}{dt}t+F_0(0),$$
for a positive constant $\delta>0$. As $1<p<p_c(n)$, it is easy to check that $n+1-(n-1)p/2>1$. Hence the following estimate is valid when $t$ is sufficiently large:
\begin{equation}\label{equation4}
F_0(t)\geq \frac{1}{2} \delta(t+R)^{n+1-(n-1)p/2}.
\end{equation}
Estimates (\ref{equation2})-(\ref{equation4}) and Lemma \ref{lemma1} with parameters
$$a\equiv n+1-(n-1)p/2,\quad\mbox{and}\quad q\equiv n(p-1)$$
imply Theorem \ref{theo1} for all exponents $p$ such that
$$(p-1)(n+1-(n-1)p/2)>n(p-1)-2,\quad\mbox{and}\quad p>1.$$
Note that the last condition on $p$ is equivalent to $p\in(1,p_c(n))$.\\

\noindent \underline{\textbf{Case $n=1$}}: In one dimensional case, the exterior domain is reduced on the semi-infinite interval $[0,\infty)$. By lemma \ref{lemma5}, we have $0<\phi_0(x)< C_2x$, then (\ref{equation5}) implies
$$
\int_0^\infty\left|u(x,t)\right|^p\phi_0(x)\,dx\geq \frac{|F_0(t)|^p}{\displaystyle \left( \int_0^{t+R}C_2x\,dx\right)^{p-1}}=C(t+R)^{-2(p-1)}|F_0(t)|^p,$$
and then by (\ref{equation3}) we get
\begin{equation}\label{equation6}
\frac{d^2}{dt^2}F_0(t)\geq k(t+R)^{-2(p-1)}|F_0(t)|^p,
\end{equation}
where $k=C^{-1}>0$. So $F_0$ satisfies the second inequality in Lemma \ref{lemma1}. To provide that $F_0$ is also verifies the first  inequality in Lemma \ref{lemma1}, we relate $\frac{d^2}{dt^2}F_0(t)$ to $F_1$ using again H\"older's inequality
\begin{eqnarray*}
\left|\int_0^\infty u(x,t)\psi_1(x,t)\,dx\right|&=&\left|\int_0^{t+R}u(x,t)[\phi_0(x)]^{1/p}[\phi_0(x)]^{-1/p}\psi_1(x,t)\,dx\right|\\
&\leq&\left( \int_0^{t+R}\left|u(x,t)[\phi_0(x)]^{1/p}\right|^p\,dx\right)^{1/p}\left( \int_0^{t+R}\left|[\phi_0(x)]^{-1/p}\psi_1(x,t)\right|^{p'}\,dx\right)^{1/{p'}}\\
&\leq&\left( \int_0^\infty\left|u(x,t)\right|^p\phi_0(x)\,dx\right)^{1/p}\left( \int_0^{t+R}[\phi_0(x)]^{-1/(p-1)}[\psi_1(x,t)]^{p'}\,dx\right)^{1/{p'}},
\end{eqnarray*}
then
\begin{eqnarray*}
|F_1(t)|^p=\left|\int_0^\infty u(x,t)\psi_1(x,t)\,dx\right|^p&\leq&\left( \int_0^\infty\left|u(x,t)\right|^p\phi_0(x)\,dx\right)\left( \int_0^{t+R}[\phi_0(x)]^{-1/(p-1)}[\psi_1(x,t)]^{p'}\,dx\right)^{p-1}.
\end{eqnarray*}
So, by using Lemmas \ref{lemma8} and \ref{lemma9}, we get
$$\int_0^\infty \left|u(x,t)\right|^p\phi_0(x)\,dx\geq \frac{|F_1(t)|^p}{\displaystyle\left( \int_0^{t+R}[\phi_0(x)]^{-1/(p-1)}[\psi_1(x,t)]^{p'}\,dx\right)^{p-1}}\geq \frac{(\varepsilon c_0)^p}{\displaystyle C^{p-1}}=L,$$
where $L>0$ is a positive contant independent of $t$. Therefore, (\ref{equation3}) implies
$$\frac{d^2}{dt^2}F_0(t)\geq L.$$
Integrate twice on $[0,t+R]$, we have
$$F_0(t)\geq \frac{1}{2}L(t+R)^{2}+\frac{dF_0(0)}{dt}t+F_0(0).$$
Hence the following estimate is valid when $t$ is sufficiently large:
\begin{equation}\label{equation7}
F_0(t)\geq \frac{1}{4}L(t+R)^{2}.
\end{equation}
Estimates (\ref{equation7}) together with (\ref{equation6}) and Lemma \ref{lemma1} with parameters
$$a\equiv 2,\quad\mbox{and}\quad q\equiv 2(p-1)$$
imply Theorem \ref{theo1} for all exponents $p$ such that
$$1<p<p_c(1)=\infty.$$

\noindent \underline{\textbf{Case $n=2$ and $p<p_c(2)$}}: As $0\not\in\overline{\Omega}$, then without loss of generality we can assume that $B_{2}(0)\cap\Omega=\emptyset$, ($B_{2}(0)$ stands for the closed ball of center 0 and radius $2$).
By lemma \ref{lemma4}, we have 
\begin{eqnarray*}
\int_{\Omega\cap\{|x|\leq t+R\}}\phi_0(x)\,dx&\leq&\int_{\{|x|\leq t+R\}\setminus B_{2}(0)}C\ln |x|\,dx\\
&=&\int_0^{2\pi}\int_{2}^{t+R}C\ln r\cdotp r\,dr\,d\theta\\
&=&\pi C\left[(t+R)^2\ln(t+R)-4\ln(2)-\frac{1}{2}((t+R)^2-4))\right]\\
&\leq&\pi C(t+R)^2\ln(t+R).
\end{eqnarray*}
Therefore,  (\ref{equation3}) and (\ref{equation5}) imply
\begin{equation}\label{equation8}
\frac{d^2}{dt^2}F_0(t)=\int_{\Omega}\left|u(x,t)\right|^p\phi_0(x)\,dx\geq \frac{|F_0(t)|^p}{\displaystyle \left( \pi C(t+R)^2\ln(t+R)\right)^{p-1}}=k[\ln(t+R)]^{-(p-1)}(t+R)^{-2(p-1)}|F_0(t)|^p
\end{equation}
where $k>0$. So $F_0$ satisfies the second inequality in Lemma \ref{lemma10}. To provide that $F_0$ is also verifies the first  inequality in Lemma \ref{lemma10}, we relate $\frac{d^2}{dt^2}F_0(t)$ to $F_1$ using again H\"older's inequality
\begin{eqnarray*}
\left|\int_{\Omega}u(x,t)\psi_1(x,t)\,dx\right|&=&\left|\int_{\Omega\cap\{|x|\leq t+R\}}u(x,t)[\phi_0(x)]^{1/p}[\phi_0(x)]^{-1/p}\psi_1(x,t)\,dx\right|\\
&\leq&\left( \int_{\Omega\cap\{|x|\leq t+R\}}\left|u(x,t)[\phi_0(x)]^{1/p}\right|^p\,dx\right)^{1/p}\left( \int_{\Omega\cap\{|x|\leq t+R\}}\left|[\phi_0(x)]^{-1/p}\psi_1(x,t)\right|^{p'}\,dx\right)^{1/{p'}}\\
&\leq&\left( \int_{\Omega}\left|u(x,t)\right|^p\phi_0(x)\,dx\right)^{1/p}\left( \int_{\Omega\cap\{|x|\leq t+R\}}[\phi_0(x)]^{-1/(p-1)}[\psi_1(x,t)]^{p'}\,dx\right)^{1/{p'}},
\end{eqnarray*}
then
\begin{eqnarray*}
|F_1(t)|^p=\left|\int_{\Omega}u(x,t)\psi_1(x,t)\,dx\right|^p&\leq&\left( \int_{\Omega}\left|u(x,t)\right|^p\phi_0(x)\,dx\right)\left( \int_{\Omega\cap\{|x|\leq t+R\}}[\phi_0(x)]^{-1/(p-1)}[\psi_1(x,t)]^{p'}\,dx\right)^{p-1}.
\end{eqnarray*}
So, by using Lemmas \ref{lemma8} and \ref{lemma9}, we get
$$\int_{\Omega}\left|u(x,t)\right|^p\phi_0(x)\,dx\geq \frac{|F_1(t)|^p}{\displaystyle\left( \int_{\Omega\cap\{|x|\leq t+R\}}[\phi_0(x)]^{-1/(p-1)}[\psi_1(x,t)]^{p'}\,dx\right)^{p-1}}\geq \frac{(\varepsilon c_0)^p}{\displaystyle \left( C(t+R)^{1-p'/2}\right)^{p-1}}=L(t+R)^{-(p/2-1)},$$
where $L>0$ is a positive contant independent of $t$. Therefore, (\ref{equation3}) implies
$$\frac{d^2}{dt^2}F_0(t)\geq L(t+R)^{-(p/2-1)}.$$
Integrate twice, we have
$$F_0(t)\geq \delta(t+R)^{3-p/2}+\frac{dF_0(0)}{dt}t+F_0(0),$$
for a positive constant $\delta>0$. As $1<p<p_c(2)$, it is easy to check that $3-p/2>1$. Hence the following estimate is valid when $t$ is sufficiently large:
\begin{equation}\label{equation4}
F_0(t)\geq \frac{1}{2} \delta(t+R)^{3-p/2}.
\end{equation}
Estimates (\ref{equation8}), (\ref{equation4}) and Lemma \ref{lemma10} with parameters
$$a\equiv 3-p/2,\quad\mbox{and}\quad q\equiv 2(p-1)$$
imply Theorem \ref{theo1} for all exponents $p$ such that
$$(p-1)(3-p/2)>2(p-1)-2,\quad\mbox{and}\quad p>1.$$

\noindent \underline{\textbf{Case $n=2$ and $p=p_c(2)$}}: As the subcritical case ($p<p_c(2)$), we have 
\begin{equation}\label{equation10}
\frac{d^2}{dt^2}F_0(t)\geq K_1[\ln(t+R)]^{-(p-1)}(t+R)^{-2(p-1)}|F_0(t)|^p,
\end{equation}
where $K_1>0$, and 
$$
\frac{d^2}{dt^2}F_0(t)=\int_{\Omega}\left|u(x,t)\right|^p\phi_0(x)\,dx\geq \frac{|F_1(t)|^p}{\displaystyle\left( \int_{\Omega\cap\{|x|\leq t+R\}}[\phi_0(x)]^{-1/(p-1)}[\psi_1(x,t)]^{p'}\,dx\right)^{p-1}}.
$$
Next, we use Lemma \ref{lemma9} and the fact that (see Lemma \ref{lemma8})
$$\int_{\Omega\cap\{|x|\leq t+R\}}[\phi_0(x)]^{-1/(p-1)}[\psi_1(x,t)]^{p'}\,dx\leq C(t+R)^{1-p'/2}(\ln(t+R))^{-1/(p-1)}.$$
we conclude that
\begin{equation}\label{equation11}
\frac{d^2}{dt^2}F_0(t)\geq \frac{\varepsilon^pc_0^p}{\displaystyle\left( C(t+R)^{1-p'/2}(\ln(t+R))^{-1/(p-1)}\right)^{p-1}}\geq L(t+R)^{-(p/2-1)}(\ln(t+R))
\end{equation}
where $L>0$ is a positive contant independent of $t$. Integrate twice, we have when $t$ is sufficiently large:
$$
F_0(t)\geq C(t+R)^{3-p/2}\ln t.
$$
As $\displaystyle \lim_{t\rightarrow\infty}\ln t=\infty$, we infer that
\begin{equation}\label{equation12}
F_0(t)\geq K_0(t+R)^{3-p/2},
\end{equation}
with $K_0>0$ being arbitrarily large when $t$ is sufficiently large.
Estimates (\ref{equation12}) together with (\ref{equation10}) and Lemma \ref{lemma11} with parameters
$$a\equiv 3-p/2,\quad\mbox{and}\quad q\equiv 2(p-1)$$
imply Theorem \ref{theo1}, since exponent $p=p_c(2)$ satisfies
$$(p-1)(3-p/2)=2(p-1)-2,\quad\mbox{and}\quad p>1.$$
\hfill$\square$



%% References

\section*{Acknowledgements}
The author would like to express sincere gratitude to Professor Ryo Ikehata for valuable discussion. 
%The referee deserves thanks for careful reading and many useful comments leading to improvement of the paper.

\bibliographystyle{elsarticle-num}
%\bibliography{bibliographie(damping2)}



\section*{References}
 \bibitem{Han} {Wei Han}, Concerning the Strauss Conjecture for the subcritical and
critical cases on the exterior domain in two space dimensions, Nonlinear Analysis ${\bf 84}$ $(2013),$ $136-145.$

\bibitem{Han2} {Wei Han}, Blow Up of Solutions to One Dimensional
Initial-Boundary Value Problems for Semilinear Wave Equations with Variable Coefficients, J. Part. Diff. Eq. ${\bf 26}$ $(2013),$ No. 2 $138-150.$ 

 \bibitem{IkehataInoue} {R. Ikehata, Yu-ki Inoue}, Global existence of weak solutions for two-dimensional semilinear wave equations with strong damping in an exterior domain, Nonlinear Analysis ${\bf 68}$ $(2008),$ $154-169$. 

 \bibitem{Jiao} {H. Jiao, Z. Zhou}, An elementary proof of the blow up for semilinear wave equation in high space dimensions, J. Differential Equations ${\bf 189}$ $(2003),$ $335-365$.

 \bibitem{Glassey} {R.T. Glassey}, Finite-time blow-up for solutions of nonlinear wave equations, Math. Z. ${\bf 177}$ $(1981),$ $323-340.$

\bibitem{Sideris}{T.C. Sideris}, Nonexistence of global solutions to semilinear wave equations in high dimensions, J. Differential Equations ${\bf 52}$ $(1984),$ $378-406.$

\bibitem{Strauss1}{W.A. Strauss}, Nonlinear Wave Equations, CBMS Reg. Conf. Ser. Math., vol. ${\bf 73}$ , AMS, Providence, RI, $1989$.

\bibitem{Strauss2}{W.A. Strauss}, Nonlinear scattering theory at low energy, J. Funct. Anal. ${\bf 41}$ $(1981)$ $110-133.$

%\bibitem{YordanovZhang} {B. Yordanov, Q.S. Zhang}, Finite time blow up for critical wave equations in high dimensions, J. Funct. Anal. ${\bf 231}$ $(2006),$ $361-374.$ 

\bibitem{ZhouHan} {Yi Zhou, Wei Han}, Blow-up of solutions to semilinear wave equations with variable
coefficients and boundary, J. Math. Anal. Appl. ${\bf 374}$ $(2011),$ $585-601.$

\end{document}